\documentclass[12pt]{article}
\usepackage{hyperref} 
\usepackage{url}
\usepackage{color} 
\usepackage{fancyvrb}
\usepackage{amsfonts}
\usepackage{amsmath}
\usepackage[affil-it,auth-sc]{authblk}
\usepackage{tabularx}
\begin{document}

\newtheorem{theorem}{Theorem}
\newtheorem{lemma}[theorem]{Lemma}
\newtheorem{proposition}[theorem]{Proposition}
\newtheorem{corollary}[theorem]{Corollary}

\newcommand{\qed}{\nobreak \ifvmode \relax \else
      \ifdim\lastskip<1.5em \hskip-\lastskip
      \hskip1.5em plus0em minus0.5em \fi \nobreak
      $\Box$\fi}

\title{\bf How to Lose with Least Probability}
\author[1]{Robert W. Chen}
\author[2]{Burton Rosenberg}
\affil[1]{Department of Mathematics, University of Miami}
\affil[2]{Department of Computer Science, University of Miami}
\maketitle

\begin{abstract}
Two players alternate tossing a biased coin where the 
probability of getting heads is $p$. The current player
is awarded $\alpha$ points for tails and $\alpha+\beta$ for heads. 
The first player reaching $n$ points wins. 
For a completely unfair coin the player going first certainly wins. 
For other coin biases, the player going first has the advantage,
but the advantage depends on the coin bias. We calculate the first 
player's advantage and the coin bias minimizing this advantage.
~

\end{abstract}

\section{Introduction}

Alice and Bob alternate tossing a biased coin where the 
probability of getting heads is $p$. Players alternate turns,
with Alice going first. The player tosses the coin and
is awarded $\alpha$ points for tails and $\alpha+\beta$ for heads. 
The first player accumulating $n$ points wins the game.
The values $\alpha, \beta $ and $n$ are positive constants.

In this paper we investigate $I(p\,|\,n,\alpha,\beta)$, 
the winning probability for Alice, the first player.
We show that $I(p\,|\,n,\alpha,\beta)$ is a polynomial in $p$ of degree 
$2m - 2$, where $m$ is the smallest positive integer greater than or equal
to $n/\alpha$.

It is easy to see that Alice has the advantage in this game since she goes first. 
As a trivial example, when $p$ is 0 or 1 then 
$I(p\,|\,n,\alpha, \beta) = 1$. We consider $p^*_n$, the value of $p$ that 
minimizes $I(p\,|\,n,\alpha, \beta)$, and hence makes the game most
favorable to Bob,
\[
I(p^*_n\,|\,n,\alpha,\beta) =\inf_{0\leq p \leq 1}I(p\,|\,n,\alpha,\beta).
\]
We provide examples of specific polynomials
$I(p\,|\,n,\alpha \beta)$ and values of $p^*_n$, generated by Mathematica.
Since these are cumbersome to calculate, and not closed form, we find
a simple form for $p^*$, the limit of $p^*_n$ as $n$ goes to infinity,
\[
p^* = \lim_{n \rightarrow \infty}  p^*_n = 1 + t - \sqrt{1 + t + t^2},
\]
where $t=\alpha/\beta$.

\section{Theorems}


\begin{theorem}[Main Theorem]
Let $\alpha, \beta>0$ and the limit $n>0$.
Then the probability that the first player wins the game is, 
\[
I(p\,|\,n,\alpha,\beta) = \frac{1}{2}\left(1+ \sum_k P(\tau_1 = k)^2\right),
\]
where the sum can range only over
$\lceil n/(\alpha + \beta) \rceil \le k \le \lceil n/\alpha \rceil$ and,
\[
P(\tau_1 = k ) =  
 {{k-1} \choose {i_{k}}}
            p^{i_{k}+1}q^{k-i_{k}-1} 
                   + \sum ^{i^*_k}_{j=i_k+1} {{k-1} \choose {j}} p^jq^{k-j-1},
\]
with $i_k=\lceil (n - k\alpha)/\beta - 1 \rceil$ and 
$i^*_k =\lceil (n - (k-1)\alpha)/\beta  - 1 \rceil$.
\end{theorem}

\noindent{\bf Proof:} 
Define stopping times 
$\tau_1$ and $\tau_2$, where $\tau_1$ is the first time that the first
player has $n$ 
or more accumulated points and $\tau_2$ is the first time that the
second player has $n$ or more accumulated points. 
Since the first player leads, 
\[
I(p\,|\,n,\alpha,\beta) = P(\tau_1 < \tau_2) +P(\tau_1 = \tau_2).
\]
Since,
\[
P(\tau_1 < \tau_2) + P(\tau_2 < \tau_1) + P(\tau_1 = \tau_2) = 1
\]
and by symmetry,
\[
P(\tau_1 < \tau_2) = P(\tau_2 < \tau_1),
\]
then,
 \[
 I(p\,|\,n,\alpha,\beta) = \bigl(1 + P(\tau_1 = \tau_2)\bigr)/2.
 \]
Note,
\begin{align*}
P(\tau_1 = \tau_2) &  = \sum_k P(\tau_1= k \,\land \,\tau_2 = k) \\
& = \sum_k P(\tau_1 = k)P(\tau_2 = k) 
= \sum_k P(\tau_1 = k)^2
\end{align*}
since $\tau_1,$ $\tau_2$ are independent and they have the 
same distribution. So we need only to compute $P(\tau_1 = k).$ 

Since unless $\lceil n/(\alpha + \beta)  \rceil$ turns are taken, player one 
cannot win, and if more than $ \lceil n/\alpha \rceil$ turns are
taken player one must have already won, we have that $P(\tau_1=k)=0$
unless $\lceil n/(\alpha + \beta) \rceil \le k \le \lceil n/\alpha \rceil$, hence the
above sums can be restricted to that range.

Consider the case where the player has not won after $k-1$ tosses,
and that $i$ heads are tossed in the first $k-1$ tosses, that the
player wins on the $k$-th toss, and that the $k$-th toss is a tail. 
This is possible for all $i$ such that,
\[
(i\beta + (k-1)\alpha < n) \;\land\; (i\beta + k\alpha \ge l) \;\land\; i\in [0,k-1].
\]
Since it does not matter the placement of the heads in the first $k-1$ tosses,
each $i$ contributes to the probability $P(\tau_1 = k)$ the amount,
\[
{{k-1} \choose i}p^iq^{k-1-i}q.
\]
Consider the similar situation, but where the $k$-th toss is a head. The constraint
on $i$ is,
\[
(i\beta+(k-1)\alpha < n ) \;\land\; ((i+1)\beta + k\alpha\ge n)\;\land\; i\in [0,k-1],
\]
and for each $i$ the contribution to $P(\tau_1 = k)$ is the amount,
\[
{{k-1} \choose i}p^iq^{k-1-i}p.
\]
Since these two cases share the range of $i$ such that,
\[
n-k\alpha \le i\beta < n-(k-1)\alpha
\]
such $i$ contribute,
\[
{{k-1} \choose i}p^iq^{k-1-i}(p+q) = {{k-1} \choose i}p^iq^{k-1-i}.
\]
The remaining $i$ must have a head on the $k$-th toss, and must satisfy,
\[
n-k\alpha-\beta \le i\beta < n-k\alpha
\]
and contribute,
\[
{{k-1} \choose i}p^{i+1}q^{k-1-i}.
\]
The only possible such $i$ is,
\[
i_k = \lceil (n-k\alpha)/\beta -1 \rceil
\]
when $i_k$ it is in the range $[0,k-1]$.
However, the binomial coefficient is zero outside
this range, so we can drop this requirement in the writing of the formula for 
$P(\tau_1 = k)$.

Returning to the case where either a heads or tails is possible on the $k$-th toss, the 
inequalities rearrange to,
\[
\lceil (n-k\alpha)/\beta \rceil \le i \le \lceil (n-(k-1)\alpha)/\beta \rceil -1 = i_k^*.
\]
Again, we can drop the requirement that such $i$ be in the range $[0,k-1]$ 
in writing the formula for $P(\tau_1 = k)$ since the
binomial coefficient is zero for $i$ outside this range.
\hfill$\Box$

\noindent{\em Note:}
It is possible that $i_k^* = i_k$, and therefore the large summation can be
empty.

\begin{corollary}
When $\lceil n/(\alpha + \beta) \rceil = \lceil n/\alpha \rceil$
then $I(p\,|\,n,\alpha,\beta) = 1$ for  all coin bias $p.$
\end{corollary}

\noindent{\bf Proof:} 
In the previous proof, the range of $k$ to sum is only the single value,
\[
k = \lceil n/(\alpha + \beta) \rceil = \lceil n/\alpha \rceil.
\]
From the calculation of $k$ we have $k\alpha\ge n$ and $k(\alpha+\beta)< n+ \alpha + \beta$.
Hence $i_k<0$ and 
 for any $i\in [0,k-1]$, $i$ heads among the first $k-1$ tosses does not win, but any
$k$ tosses does win. Therefore the sum for $P(\tau_1 = k )$ reduces to, 
\begin{align*}
P(\tau_1 = k ) & =  
 {{k-1} \choose {i_{k}}}
            p^{i_{k}+1}q^{k-i_{k}-1} 
                   + \sum ^{i^*_k}_{j=i_k+1} {{k-1} \choose {j}} p^jq^{k-j-1}\\
       &= 0 +  \sum ^{k-1}_{j=0} {{k-1} \choose {j}} p^jq^{k-j-1} = 1.
\end{align*}
\hfill$\Box$

\begin{corollary}
$ I(p\,|\,n,\alpha,\beta)$ is either a polynomial in $p$ of degree $2\lceil n/\alpha \rceil-2$,
or is 1, in the case $\lceil n/(\alpha + \beta) \rceil = \lceil n/\alpha \rceil$.
\end{corollary}

\noindent{\bf Proof:} 
The $P(\tau_1 = k )$ are polynomials  in $p$ of highest degree $\lceil n/\alpha \rceil-1$,
since $i_k<i^*_k$, and $i^*_k$ is the maximum number of heads among the first $k-1$
tosses yet the player does not win, and this is bounded by $\lceil n/\alpha\rceil-1$ which
is the maximum number of toss of all tails such that the player does not win.
\hfill$\Box$

\begin{theorem}[Optimal coin bias]
Assumptions as in the previous theorem, let $p^*_n$ be the coin bias that
minimizes the probability of the first player winning. Then, 
\[
p^* = \lim_{n \rightarrow \infty}  p^*_n = 1 + t - \sqrt{1 + t + t^2},
\]
where $t=\alpha/\beta$.
\end{theorem}

\noindent{\bf Proof:} 
In order to find the limit $p^*$ of $p^*_n$ as $n 
\rightarrow 
\infty,$ let $X_1, X_2, \ldots $ be i.i.d. random variables such that $P(X_i = 1) 
= p$ and $P(X_i = 0) = 1 - p.$ For each $k,$ let 
$S_k = X_1 + \ldots +X_k$ and  $P(\tau _1 > k) = P(k\alpha + \beta 
S_k < n).$ 

For the positive integer $k$ choose $x$ such that 
$k = n/(\alpha + p\beta ) + \sqrt{n}x$. 
Then,
\begin{align*}
k\alpha + \beta S_k 
&\approx k\alpha + \beta \left( kp +  \sqrt{kp(1-p)}Z\right) \\
&\approx  n + (\alpha + \beta p) \sqrt{n}x + \beta \sqrt{\frac{np(1-p)}{\alpha + \beta p}}Z.
\end{align*}
where $Z$ is  the standard normal random variable. 
Therefore,
\[
P(k\alpha + \beta S_k < n) \approx 
P(Z < - \sigma x),
\]
 where 
\[
\sigma ^2 = (\alpha + \beta p)^3/(\beta ^2 
p(1-p)).
\]
Since,
\[
P(Z < - \sigma x) \approx \int _{-\infty}^{-\sigma x}1/\sqrt{2
\pi }e^{-u^2/2}du,
\]
we have,
\[
P(\tau _1 = \tau _2) = \sum _kP^2(\tau _1 = k) \approx 
\sigma ^2 \int _{-\infty}^{\infty}e^{-\sigma ^2 u^2}du.
\]
It is easy to see 
that to minimize $I(p\,|\,n,\alpha ,\beta )$ is the same to minimize $P(\tau _1 
= \tau _2).$ Since the integral on the right hand side of the above equation
is increasing in $\sigma ^2$, to minimize $I(p\,|\,n,\alpha ,\beta )$ we  
minimize $\sigma ^2$ with respect to $p.$ Solving,
\[d(\ln(\sigma ^2))/dp = 3\beta /(\alpha + \beta p) - 1/p + 1/(1 - p) =0
\]
so,
\[p^* = \left(\beta + 
\alpha - \sqrt{\beta ^2 + \beta \alpha + \alpha ^2}\right)/\beta = 1 + t - 
\sqrt{1 + t + t^2},
\]
where $0 < t = \alpha /\beta < \infty.$ 
\hfill$\Box$

~\\
\noindent
{\em Remarks:}
Since $\sqrt{1 + t + t^2} > 1/2 + t$ for all $0 < t < \infty,$ 
therefore $p^* < 1/2.$ For 
example, if $\alpha = \beta > 0,$ $p^* = 2 - \sqrt{3} \approx 0.267949192,$ if 
$\alpha /\beta = 2,$ $p^* = 3 - \sqrt{7} \approx 0.354248688,$ and if $\alpha /
\beta = 0.5,$ $p^* = (3 - \sqrt{7})/2 \approx  0.177124344.$ When $t \rightarrow 
\infty,$ $p^* \rightarrow 0.5$ and when $t \rightarrow 0,$ $p^* 
\rightarrow 0.$

\section{Computer experiments}

\subsection{Mathematica Code}
In Figure~\ref{code:advantage} we give Mathematica code which
generates $I(p\,|\,n,\alpha ,\beta ).$

In Figure~\ref{code:minadvantage} we give Mathematica code which
generates $I(p^*_n\,|\,n,\alpha ,\beta ),$ the coin bias that minimizes player one's
advantage.

\subsection{Advantage polynomials}
 

Table~\ref{Advantage:Alpha1:Beta1} through Table~\ref{Advantage:Alpha3:Beta2} give
example  $I(p\,|\,n,\alpha ,\beta)$ for various  values of $n,$ $\alpha $, and $\beta $.

\subsection{Example minimized advantages}

There is an online computer program at:

\url{http://epaper-live.appspot.com/2011-CGZ} 

\noindent using a randomized  simulation to estimate $I(p^*_n\,|\,n,\alpha,\beta).$

The following are some values of $I(p^*_n\,|\,n,\alpha ,\beta )$ and 
$I(p^*\,|\,n,\alpha ,\beta )$ for $m = 5,\;10,\;15$ and various values of $\alpha $ and
$\beta .$

Notice that when $m$ is moderate or large, it is very cumbersome to write down 
the polynomial of $I(p\,|\,n,\alpha, \beta).$  Since $I(p^*\,|\,n,\alpha,\beta)$ is a 
good approximation for $I(p^*_n\,|\,n,\alpha,\beta),$ we may use $p^*$ for the 
probability for heads and use simulation to estimate $I(p^*_n\,|\,n,\alpha,\beta).$

\begin{figure}\begin{Verbatim}[label=Mathematica Code, tabsize=3, fontfamily=tt, numbers=left, numbersep=4pt, frame=single, fontsize=\footnotesize]
PlayerOne[n_, alpha_, beta_] := Module[{l, m, k, expr},
	l = Ceiling[n/(alpha + beta)];
	m = Ceiling[n/alpha];
	For[k = l, k <= m, k++, Module[{h, i, j, k = k},
		h = Max[0, Ceiling[(n - k*alpha)/beta - 1]];
		i = Max[0, Ceiling[(n - k*alpha + alpha)/beta - 1]];
		If[h == i, p[k, p_] := 
		    Binomial[k - 1, h]*p^(h + 1)*(1 - p)^(k - h - 1)];
		If[h < i, p[k, p_] := 
		  Binomial[k-1,h]*p^(h + 1)*(1 - p)^(k - h - 1) 
		  + Sum[Binomial[k-1,j]*p^j*(1 - p)^(k - j - 1), {j, h + 1, i}]];
		If[k == m, p[m, p_] := (1 - p)^(m - 1) 
		  + Sum[Binomial[m - 1, j]*p^j*(1 - p)^(m - j - 1), {j, 1, i}]];
		]];
	expr = Expand[(1/2)*(1 + Sum[p[k, p]^2, {k, l, m}])];
	If[l == m, expr = 1;];
	Return[expr]
];
\end{Verbatim}
\caption{Polynomial expression for the advantage of player one.}
\label{code:advantage}
\end{figure}

\begin{figure}\begin{Verbatim}[label=Mathematica Code, tabsize=3, fontfamily=tt, numbers=left, numbersep=4pt, frame=single, fontsize=\footnotesize]
PlayerOneMin[n_, alpha_, beta_] := 
   N[Minimize[Evaluate[PlayerOne[n, alpha , beta ]], p]]
\end{Verbatim}
\caption{Minimizing the advantage of player one.}
\label{code:minadvantage}
\end{figure}

\newcolumntype{R}{>{\raggedright\arraybackslash}X}
\begin{table}
\begin{tabularx}{\textwidth}{|r|R|}
\hline
$n$& 
$I(p\,|\,n,\alpha=1,\beta=1)$ \\ \hline\hline
1&
$1$ \\ \hline
2&
$1 - p + p^2$ \\ \hline
3&
$1 -2p + 5p^2 -4p^3 + p^4$ \\ \hline
4&
$1 - 3p + 12p^2 - 22p^3 + 19p^4 - 7p^5 + p^6$ \\ \hline
5&
$1 - 4p + 22p^2 - 64p^3 + 104p^4 - 92p^5 + 43p^6 - 10p^7 + p^8$\\ \hline
6&
$1 - 5p + 35p^2 - 140p^3 + 341p^4 - 508p^5 
+ 459p^6 - 247p^7 + 77p^8 - 13p^9 + p^{10}$\\ \hline
7&
$1 - 6p + 51p^2 -260p^3 + 850p^4 - 1,816p^5 
+ 2,548p^6 - 2,336p^7 + 1,385p^8  - 522p^9 + 121p^{10} - 16p^{11} + p^{12}$\\ \hline
8&
$1 - 7p + 70p^2 - 434p^3 + 1,786p^4 - 5,011p^5 + 9,709p^6 
- 13,030p^7 + 12,079p^8 -7,683p^9 + 3,316p^{10} - 952p^{11} 
+ 175p^{12} - 19p^{13} + p^{14}$\\ \hline
9&
$1 - 8p + 92p^2 - 672p^3 + 3,339p^4 -11,648p^5 + 29,037p^6 
- 52,154p^7 + 67,644p^8 -  63,248p^9 + 42,440p^{10} - 20,280p^{11} 
+ 6,812p^{12} - 1,572p^{13} + 239p^{14} - 22p^{15} + p^{16}$\\ \hline
10&
$1 - 9p + 117p^2 - 984p^3 + 5,734p^4 - 23,968p^5 
+ 73,381p^6 - 166,489p^7 + 281,524p^8 - 355,393p^9 
+ 334,585p^{10} - 234,160p^{11} + 121,147p^{12} - 45,916p^{13} 
+ 12,559p^{14} - 2,417p^{15} + 313p^{16} - 25p^{17} + p^{18}$\\ \hline
11&
$1 - 10p + 145p^2 - 1,380p^3 + 9,231p^4 - 45,024p^5 
+ 163,864p^6 - 451,312p^7+ 948.352p^8 - 1,526,710p^9 + 1,885,453p^{10} 
- 1,785,028p^{11} + 1,292,464p^{12} - 712,744p^{13} + 297,382p^{14} 
- 92,900p^{15} + 21,369p^{16} - 3,522p^{17} + 397p^{18} - 28p^{19} + p^{20}$\\ \hline
12&
$1 - 11p + 176p^2 - 1,870p^3 + 14,125p^4 - 78,807p^5 
+ 332,865p^6 - 1,081,308p^7+ 2,728,525p^8 - 5,379,992p^9 
+ 8,314,959p^{10} - 10,083,869p^{11} + 9,591,305p^{12} - 7,142,250p^{13} 
+ 4,150,664p^{14} - 1,873,073p^{15} + 651,365p^{16} - 172,523p^{17} 
+ 34,180p^{18} - 4,922p^{19} + 491p^{20} - 31p^{21} + p^{22}$ \\  \hline
\end{tabularx}
\caption{$\alpha=1$, $\beta=1$}
\label{Advantage:Alpha1:Beta1}
\end{table}

{\newcolumntype{R}{>{\raggedright\arraybackslash}X}
\begin{table}
\begin{tabularx}{\textwidth}{|r|R|}
\hline
$n$& 
$I(p\,|\,n,\alpha=1,\beta=2)$ \\ \hline\hline 
1&$1$ \\ \hline
2& $1 - p + p^2$ \\ \hline
3& $1 - 2p + 4p^2 - 3p^3 + p^4$ \\ \hline
4&$1 - 3p + 10p^2 - 14p^3 + 11p^4 - 5p^5 + p^6$ \\ \hline
5&$1 - 4p + 19p^2 - 43p^3 + 54p^4 - 42p^5 + 22p^6 
- 7p^7 +p^8$ \\ \hline
6&$1 - 5p + 31p^2 - 100p^3 + 186p^4 - 216p^5 + 
169p^6 - 94p^7 + 37p^8 - 9p^9 + p^{10}$ \\ \hline
7&$1 - 6p + 46p^2 - 195p^3 + 497p^4 - 808p^5 + 
891p^6 - 700p^7 + 407p^8 - 178p^9 + 56p^{10} - 11p^{11} + p^{12}$ \\ \hline
8&$1 - 7p + 64p^2 - 338p^3 + 1,112p^4 - 2,393p^5 +
3,538p^6 - 3,755p^7 + 2,963p^8 - 1,785p^9 + 836p^{10} - 302p^{11} + 
79p^{12} - 13p^{13} + p^{14}$ \\ \hline
9&$1 - 8p + 85p^2 - 539p^3 + 2,191p^4 - 5,966p^5 +
11,335p^6 - 15,606p^7 + 16,088p^8 - 12,744p^9 + 7,911p^{10} - 3,905p^{11} 
+ 1,540p^{12} - 474p^{13} + 106p^{14} - 15p^{15} + p^{16}$ \\  \hline
10&$1 - 9p + 109p^2 - 808p^3 + 3,929p^4 - 13,068p^5 
+ 30,811p^6 - 53,204p^7 + 69,302p^8 - 69,820p^9 + 55,500p^{10} - 35,338p^{11}
+ 18,228p^{12} - 7,670p^{13} + 2,619p^{14} - 702p^{15} + 137p^{16} - 
17p^{17} + p^{18}$ \\ \hline
11&$1 - 10p + 136p^2 - 1,155p^3 + 6,556p^4 - 25,912p^5
+ 73,689p^6 - 155,186p^7 + 248,468p^8 - 309,592p^9 + 306,164p^{10} - 244,118p^{11}
+ 158,834p^{12} - 85,084p^{13} + 37,759p^{14} - 13,898p^{15} + 4,189p^{16} -
994p^{17} + 172p^{18} - 19p^{19} + p^{20}$ \\ \hline
12&$1 - 11p + 166p^2 - 1,590p^3 + 10,337p^4 - 47,509p^5
+ 159,238p^6 - 399,557p^7 + 768,438p^8 - 1,157,198p^9 + 1,390,348p^{10} - 1,354,014p^{11}
+ 1,082,537p^{12} - 717,534p^{13} + 397,153p^{14} - 184,505p^{15} + 72,133p^{16} -
23,647p^{17} + 6,382p^{18} - 1,358p^{19} + 211p^{20} - 21p^{21} + p^{22}$ \\ \hline
\end{tabularx}
\caption{$\alpha=1$, $\beta=2$}
\label{Advantage:Alpha1:Beta2}
\end{table}

\newcolumntype{R}{>{\raggedright\arraybackslash}X}
\begin{table}
\begin{tabularx}{\textwidth}{|r|R|}
\hline
$n$& 
$I(p\,|\,n,\alpha=2,\beta=1)$ \\ \hline\hline 
2&$1$ \\ \hline
4&$ 1$ \\ \hline
6&$1 - p^2 + p^4$ \\ \hline
8&$1 - 3p^2 + 2p^3 + 9p^4 - 12p^5 + 4p^6$ \\ \hline
10&$1 - 6p^2 + 8p^3 + 33p^4 - 96p^5 + 100p^6 - 
48p^7 + 9p^8$ \\ \hline
12&$1 - 10p^2 + 20p^3 + 85p^4 - 396p^5 +
690p^6 - 660p^7 + 371p^8 - 116p^9 + 16p^{10}$ \\ \hline
14&$1 - 15p^2 + 40p^3 + 180p^4 - 1,176p^5 +
2,870p^6 - 4,060p^7 + 3,735p^8 - 2,300p^9 + 921p^{10} 
- 220p^{11} + 25p^{12}$ \\ \hline
16&$1 - 21p^2 + 70p^3 + 336p^4 - 2,856p^5 + 8,960p^6 - 16,668p^7 + 
21,015p^8 - 19,032p^9 + 12,531p^{10} - 5,850p^{11} + 1,839p^{12} - 360p^{13} + 36p^{14}$ \\ \hline
18&$1 - 28p^2 + 112p^3 + 574p^4 - 6,048p^5 + 23,184p^6 - 53,264p^7 + 
84,616p^8 - 100,128p^9 + 91,385p^{10} - 64,204p^{11} + 33,678p^{12} - 12,608p^{13} + 
3,214p^{14} - 532p^{15} + 49p^{16}$ \\ \hline
20&$1 - 36p^2 + 168p^3 + 918p^4 - 11,592p^5 + 52,500p^6 - 143,256p^7 + 
273,133p^8 - 396,296p^9 + 461,258p^{10} - 437,572p^{11} + 331,905p^{12} - 193,564p^{13} + 
82,979p^{14} - 25,046p^{15} + 5,165p^{16} - 728p^{17} + 64p^{18}$ \\ \hline
22&$1 - 45p^2 + 240p^3 + 1,395p^4 - 20,592p^5 + 107,580p^6 - 339,480p^7 +
752,661p^8 - 1,287,280p^9 + 1,818,666p^{10} - 2,192,040p^{11} + 2,229,076p^{12} - 
1,842,576p^{13} + 1,187,376p^{14} - 573,496p^{15} + 199,623p^{16} - 48,144p^{17} + 7,891p^{18} 
 - 936p^{19} + 81p^{20}$ \\ \hline
 24&$1 - 55p^2 + 330p^3 + 2,035p^4 - 34,452p^5 + 203,940p^6 - 729,960p^7 +
1,840,785p^8 - 3,613,360p^9 + 5,997,570p^{10} - 8,828,788p^{11} + 11,452,248p^{12} -
12,571,976p^{13} + 11,200,776p^{14} - 7,840,824p^{15} + 4,195,879p^{16} - 1,666,048p^{17} + 
472,645p^{18} - 91,446p^{19} + 11,741p^{20} - 1,140p^{21} + 100p^{22}$ \\ \hline
\end{tabularx}
\caption{$\alpha=2$, $\beta=1$}
\label{Advantage:Alpha2:Beta1}
\end{table}

\newcolumntype{R}{>{\raggedright\arraybackslash}X}
\begin{table}
\begin{tabularx}{\textwidth}{|r|R|}
\hline
$n$& 
$I(p\,|\,n,\alpha=2,\beta=3)$ \\ \hline\hline 
2&$1$ \\   \hline
4& $1 - p + p^2$ \\ \hline
6&$1 - 2p + 5p^2 - 4p^3 +  p^4$ \\ \hline
8&$1 - 3p + 12p^2 - 22p^3 + 19p^4 - 7p^5 + p^6$ \\ \hline
10&$1 - 4p + 22p^2 - 64p^3 + 102p^4 - 90p^5 + 43p^6   
- 10p^7 + p^8$ \\ \hline
12&$1 - 5p + 35p^2 - 140p^3 + 332p^4 - 478p^5 +
429p^6 - 238p^7 + 77p^8 - 13p^9 + p^{10}$ \\ \hline
14&$1 - 6p + 51p^2 - 260p^3 + 826p^4 - 
1,678p^5 + 2,251p^6 - 2,039p^7 + 1,247p^8 - 498p^9 + 
121p^{10} - 16P^{11} + p^{12}$ \\ \hline
16&$1 - 7p + 70p^2 - 434p^3 + 1,736p^4 -
4,601p^5 + 8,335p^6 - 10,604p^7 + 9,653p^8 - 6,309p^9 +
2,906p^{10} - 902P^{11} + 175p^{12} - 19p^{13} + p^{14}$ \\ \hline
18&$1 - 8p + 92p^2 - 672p^3 + 3,249p^4 -
10,688p^5 + 24,647p^6 - 40,904p^7 + 49,930p^8 - 45,534p^9 +
31,190p^{10} - 15,890P^{11} + 5,852p^{12} - 1,482p^{13} + 
239p^{14} - 22p^{15} + p^{16}$ \\ \hline
20&$1 - 9p +117p^2 - 984p^3 + 5,587p^4 
- 22,036p^5 + 62,161p^6 - 128,594p^7 + 199,200p^8 - 
235,130p^9 + 214,326p^{10} - 151,840P^{11} + 83,252p^{12} - 
34,696p^{13} + 10,727p^{14} - 2,270p^{15} + 313p^{16} 
-25p^{17} + p^{18}$ \\ \hline
22&$1 - 10p + 145p^2 - 1,380p^3 + 
9,007p^4 - 41,524p^5 + 139,189p^6 - 347,593p^7 + 659462p^8 
- 966,212p^9 + 1,108,918p^{10} - 1,008,598P^{11} + 732,121p^{12} 
- 423,904p^{13} + 193,663p^{14} - 68,225p^{15} + 17,869p^{16} 
- 3,298p^{17} + 397p^{18} - 28p^{19} + p^{20}$ \\ \hline
24&$1 - 11p + 176p^2 - 1,870p^3 + 
13,801p^4 - 72,939p^5 + 284,173p^6 - 836,056p^7 + 1,892,572p^8    
- 3,346,345p^9 + 4,682,293p^{10} - 5,246,194P^{11} + 
4,755,285p^{12} - 3,512,514p^{13} + 2,118,592p^{14} - 
1,037,420p^{15} + 406,113p^{16} - 123,831p^{17} + 28,312p^{18} 
- 4,598p^{19} + 491p^{20} - 31p^{21} + p^{22}$ \\ \hline
\end{tabularx}
\caption{$\alpha=2$, $\beta=3$}
\label{Advantage:Alpha2:Beta3}
\end{table}

\newcolumntype{R}{>{\raggedright\arraybackslash}X}
\begin{table}
\begin{tabularx}{\textwidth}{|r|R|}
\hline
$n$& 
$I(p\,|\,n,\alpha=3,\beta=2)$ \\ \hline\hline 
3&$1$ \\   \hline
6&$1$ \\   \hline
9&$1 - p^2 + p^4$ \\ \hline
12&$1 - 3p^2 + 2p^3 + 9p^4 - 12p^5 + 4p^6$ \\ \hline
15&$1 - 6p^2 + 8p^3 + 33p^4 - 102p^5 + 109p^6 - 51p^7 + 9p^8$ \\     \hline
18&$1 - 10p^2 + 20p^3 + 85p^4 - 436p^5 + 826p^6 - 824p^7 + 455p^8 - 
132p^9 + 16p^{10}$ \\ \hline
21&$1 - 15p^2 + 40p^3 + 180p^4 - 1,326p^5 + 3,670p^6 - 5,760p^7 + 5,595p^8 
- 3,420p^9 + 1,281p^{10} - 270p^{11} + 25p^{12}$ \\ \hline
24&$1 - 21p^2 + 70p^3 + 336p^4 - 3,276p^5 + 12,020p^6 - 26,028p^7 + 
36,835p^8 - 35,347p^9 + 23,176p^{10} - 10,220p^{11} + 2,899p^{12} - 480p^{13} + 36p^{14}$ \\ \hline
27&$1 - 28p^2 + 112p^3 + 574p^4 - 7,028p^5 + 32,249p^6 - 89,524p^7 +
167,706p^8 - 221,823p^9 + 211,485p^{10} - 146,003p^{11} + 72,252p^{12} - 24,943p^{13} + 
5,699p^{14} - 777p^{15} + 49p^{16}$ \\ \hline
30&$1 - 36p^2 + 168p^3 + 918p^4 - 13,608p^5 + 75,124p^6 - 255,144p^7 +
597,093p^8 - 1,011,526p^9 + 1,274,175p^{10} - 1,210,029p^{11} + 869,107p^{12} - 468,840p^{13} + 
186,558p^{14} - 52,997p^{15} + 10,149p^{16} - 1,176p^{17} + 
64p^{18}$ \\ \hline
33&$1 - 45p^2 + 240p^3 + 1,395p^4 - 24,372p^5 + 157,476p^6 - 633,564p^7 +
1,782,165p^8 - 3,689,932p^9 + 5,794,834p^{10} - 7,033,324p^{11} + 6,664,633p^{12} - 
4,941,244p^{13} + 2,850,905p^{14} - 1,263,598p^{15} + 420,991p^{16} - 101,764p^{17} + 
16,795p^{18} - 1,692p^{19} + 81p^{20}$ \\ \hline
36&$1 - 55p^2 + 330p^3 + 2,035p^4 - 41,052p^5 
+ 304,140p^6 - 1,415,760p^7 + 4,657,305p^8 - 11,410,660p^9 + 
21,494,814p^{10} - 31,814,548p^{11} + 37,533,952p^{12} - 
35,564,264p^{13} + 27,090,180p^{14} - 16,505,584p^{15} + 
7,962,291p^{16} - 2,994,900p^{17} + 858,625p^{18} - 180,870p^{19} 
+ 26,261p^{20} - 2,340p^{21} + 100p^{22}$ \\ \hline
\end{tabularx}
\caption{$\alpha=3$, $\beta=2$}
\label{Advantage:Alpha3:Beta2}
\end{table}

\newcolumntype{R}{>{\raggedright\arraybackslash}X}
\begin{table}
\begin{tabularx}{\textwidth}{|X|X|X|X|X|}
\hline
$n$& $\alpha$ & $\beta$ &$I(p^*_n\,|\,n,\alpha ,\beta )$ & $I(p^*\,|\,n,\alpha ,\beta ) $
\\ \hline\hline
 5&1&1& 0.700& 0.700\\ \hline
10&1&1&0.643&0.643\\  \hline
15&1&1&0.617&0.617\\  \hline
5&1&2&0.689&0.695\\ \hline
10&1&2&0.635&0.641\\ \hline
15&1&2&0.610&0.616\\ \hline
5&1&3&0.660&0.663\\ \hline
10&1&3&0.624&0.630\\ \hline
15&1&3&0.604&0.609\\ \hline
10&2&1&0.750&0.753\\ \hline
20&2&1&0.714&0.718\\ \hline
30&2&1&0.683&0.686\\ \hline
10&2&3&0.669&0.701\\ \hline
20&2&3&0.714&0.719\\ \hline
30&2,&3&0.570&0.612\\ \hline
15&3&2&0.716&0.752\\ \hline
30&3&2& 0.665&0.676\\ \hline
45&3&2&0.649&0.666\\ \hline
\end{tabularx}
\caption{Sample approximate and exact advantage values}
\label{Table:approxadvantage}
\end{table}

\begin{figure}
\begin{Verbatim}[label=Python Code, tabsize=3, fontfamily=tt, numbers=left, numbersep=4pt, frame=single, fontsize=\footnotesize]
import random
import math
import os

HEADS = 1
TAILS = -1
MAX_TRIALS = 10000

def usage():
    print "Usage: game -a alpha - b beta -n limit - t trials"

def myrand(p) :
    r = random.random()
    if r < p :
        return HEADS
    return TAILS

def f(alpha, beta, n, p) :
    na = nb = 0
    turns = 0
    while 1 :
        turns += 1
        na += alpha
        a_toss = myrand(p)
        if a_toss == HEADS :
            na += beta
        if na >= n :
            return turns
        nb += alpha
        b_toss = myrand(p)
        if b_toss == HEADS :
            nb += beta
        if nb >= n :
            return -turns  
    print "never gets here"
    return
\end{Verbatim}
\caption{Python program for minimized advantage, online at URL in text.}
\label{code:online}
\end{figure}

\end{document}